\documentclass{article}
\usepackage{amssymb}
\usepackage{amscd}
\usepackage[utf8]{inputenc}
\usepackage{amsfonts}
\usepackage{amstext}
\usepackage{color}
\newtheorem{theorem}{Theorem}
\newtheorem{lemma}{Lemma}

 \newtheorem{proposition}{Proposition}
\newtheorem{corollary}{Corollary} \newcommand{\La}{\Lambda}
\newcommand{\R}{{\mathbb R}}  \newcommand{\Z}{{\mathbb Z}} \newcommand{\N}{{\mathbb N}}
\newcommand{\h}{h_\varepsilon}\newcommand{\J}{\tilde{J}}
\newcommand{\Cc}{{\mathbb C}}
\providecommand{\keywords}[1]
{
  \small	
  \textbf{\textit{Keywords: }} #1
}

\begin{document}
\author{Alexander~Ulanovskii and Ilya~Zlotnikov }

\title{Reconstruction  of Bandlimited Functions from Space--Time Samples }
\maketitle
\begin{abstract}
   For a wide family of even kernels  $\{\varphi_u, u\in I\}$, we  describe  discrete sets $\La$ such that  every bandlimited signal $f$ can be  reconstructed from the space-time samples $\{(f\ast\varphi_u)(\lambda), \lambda\in\La, u\in I\}$.
\end{abstract}
\providecommand{\keywords}[1]
{
  {\small	
  \textbf{\textit{Keywords: }} #1}
}

\date{}\maketitle
\section{Introduction} The {\it classical sampling problem} asks  when a continuous signal (function) $f$  can be reconstructed from its  discrete samples $f(\lambda),\lambda\in\La$. In the {\it dynamical sampling problem}, the set of space samples is replaced by a set of space-time samples (see e.g. \cite{A1}, \cite{A2}, \cite{g}, \cite{MMM} and references therein). An interesting case is the  problem of reconstruction of a bandlimited signal $f$ from the space-time  samples of its states $f\ast\varphi_u$ resulting from the convolution with a kernel $\varphi_u$. An important example (see \cite{g} and \cite{lv}) is the Gaussian kernel $\varphi_u(x)=\exp(-u x^2)$, which arises from  the diffusion process. More generally,  the kernel  \begin{equation}\label{ker}\varphi_u(x)=\exp(-u|x|^{\alpha}),\quad \alpha>0,\end{equation}  arises form the fractional diffusion equation.

Denote by $PW_\sigma$ the Paley--Wiener space
$$
PW_\sigma:=\{f\in L^2(\R): \mbox{supp}(\hat f)\subseteq[-\sigma,\sigma]\},
$$ where $\hat f$ denotes the Fourier transform
$$
\hat f(t)=\int\limits_\R e^{- itx}f(x)\,dx.
$$

A set $\La\subset\R$ is called uniformly discrete (u.d.) if
\begin{equation}\label{sep}
\delta(\La):=\inf_{\lambda,\lambda'\in\La,\lambda\ne\lambda'}|\lambda-\lambda'|>0.
\end{equation}

The following problem  is considered in \cite{g}: Given a u.d. set $\La\subset\R$ and  a kernel  $\{\varphi_u,u\in I\}$, where $I$  is an interval.
What are the conditions that allow one to recover a function $f\in PW_\sigma$ in a stable way from the data set
\begin{equation}\label{data}
\{(f\ast\varphi_u)(\lambda):\lambda\in\La, u\in I\}?
\end{equation}

 In what follows, we denote by $\Phi_u$ the Fourier transform of $\varphi_u$ and assume that the functions $\varphi_u(x)$ and $\Phi_u(t)$
  are continuous functions of $(x,u)$ and $(t,u)$, respectively.

  It is remarked in \cite{g}, that the property of stable recovery formulated above is equivalent to the existence of two constants
$A,B$ such that
\begin{equation}\label{ss}
A\|f\|^2_2\leq \int\limits_I\sum_{\lambda\in\La}|(f\ast\varphi_u)(\lambda)|^2\,du\leq B\|f\|_2^2,\quad \forall f\in PW_\sigma.
\end{equation}

It often happens in the sampling theory that  inequalities similar to the one in the right hand-side of (\ref{ss}) are not difficult to check.
  It is also the case here,  it suffices to assume the uniform boundedness of the $L^1(\R)$-norms $\|\varphi_u\|_1$:
     \begin{proposition}\label{p1}Assume
   \begin{equation}\label{fb}
   \sup_{u\in I}\|\varphi_u\|_1<\infty.
   \end{equation}
Then for every $\sigma>0$ and every u.d. set $\La$ there is a constant $B$ such that
$$
 \int\limits_{I}\sum_{\lambda\in\La}|(f\ast\varphi_u)(\lambda)|^2\,du\leq B \|f\|_2^2, \quad \forall f\in PW_\sigma.
$$\end{proposition}

We present a simple proof in Section 3.

Hence, the main difficulty lies in proving the left hand-side inequality.

Recall that the classical  Shannon sampling theorem states that every ${f \in PW_\sigma}$ admits a stable recovery from the uniform space samples $f(k/a), k\in\Z$, if and only if  $a\geq \sigma/\pi$. The critical value $a=\sigma/\pi$ is called the  Nyquist rate.
Since the  space-time samples (\ref{data}) produce ``more information'' compared to the space samples, one may  expect that every $f\in PW_\sigma$  can be recovered from  the space-time  uniform samples  at sub-Nyquist spatial density.  However, it is not the case, as shown  in \cite{lv} for the convolution with the Gaussian kernel.
On the other hand, it is proved in  \cite{g} that uniform dynamical samples at sub-Nyquist spatial rate allow one to stably reconstruct the Fourier transform $\hat f$ away from certain, explicitly described blind spots.

It is well-known that the nonuniform sampling is sometimes more efficient than  the uniform one. For example, this is so for the universal sampling, see  e.g. \cite{ou1}, Lecture 6. It is also the case for the problem above: For a wide class of even kernels, we  show that data (\ref{data}) always allows stable reconstruction, provided $\La$ is any relatively dense set ``different'' from an arithmetic progression.

To state  precisely our main result, we need the following definition:
Given a u.d. set $\La$, the collection of sets $W(\La)$  is defined as all weak limits of the translates  $\La-x_k$, where $x_k$ is any bounded or unbounded sequence of real numbers (for the definition of weak limit see e.g. Lecture~3.4.1 in \cite{ou1}).

Consider the following condition:

\medskip ($\alpha$)  $W(\La)$ does not contain the empty set, and no element $\La^\ast\in W(\La)$ lies in an arithmetic progression.

\medskip
The first property in ($\alpha$) means that $\La$ is relatively dense, i.e. there exists $r>0$ such that every interval $(x,x+r)$ contains at least one point of $\La.$ It follows that every element $\La^\ast\in W(\La)$ is also a relatively dense set.

The second condition in ($\alpha$) means that no $\La^\ast\in W(\La)$ is a subset of $b+(1/a)\Z$, for some $a>0$ and $b\in\R.$


     Let us now define a collection of kernels $\mathcal{C}$: A  kernel $\{\varphi_u, u\in I\}$, where $I$ is an interval, belongs to $\mathcal{C}$ if it satisfies  the following five conditions:

  \medskip

  ($\beta$) 
  There is a constant $C$ such that
  \begin{equation}\label{u1}
  \sup_{u\in I}|\varphi_u(x)|\leq\frac{C}{1+x^4},\quad x\in\R;
  \end{equation}

  ($\gamma$)  There is a constant $C$ such that
   \begin{equation}\label{u2}
   \|\varphi_{u'}-\varphi_u\|_1\leq C|u-u'|, \quad u,u'\in I;
  \end{equation}

($\zeta$)  Every $\varphi_u$ is real and even: $\varphi_u(x)\in\R, \varphi_u(-x)=\varphi_u(x),x\in\R, u\in I$;


 ($\eta$) $\sup\limits_{u \in I}|\Phi_u(t)|>0$ for every $t\in\R$;

   ($\theta$) For every $w\in\Cc$ and every $\sigma>0$, the family $\{\Phi_u''(t)+w\Phi_u(t), u\in I\}$ forms a complete set in $L^2(0,\sigma)$.

  \medskip

Clearly, condition (\ref{u1}) implies that the derivatives $\Phi_u''(t), u\in I,$ are continuous and uniformly bounded.
  Condition ($\zeta$) implies that  the functions $\Phi_u$ are real and even.

  One may easily check that $\mathcal{C}$ contains the kernels defined in  (\ref{ker}), where $I=(a,b)$ is any interval such that  $ 0<a<b<\infty$.

 Our main result is as follows:

\begin{theorem}\label{t1} Given a u.d. set  $\La\subset\R$ and a kernel $\{\varphi_u,u\in I\}\in \mathcal{C}.$ The following conditions are equivalent:

{\rm(a)} The left inequality in  {\rm(\ref{ss})} is true for every $\sigma>0$ and some $A=A(\sigma)$;

{\rm(b)} $\La$ satisfies condition {\rm($\alpha$).}
\end{theorem}

\section{Space--Time Sampling in Bernstein Spaces}

The aim of this section is to prove a variant of Theorem 1 for the Bernstein space $B_\sigma$.

It is well-know that every function $f\in PW_\sigma$ admits an analytic continuation to the complex plane and satisfies
\begin{equation}\label{f}
|f(x+iy)|\leq Ce^{\sigma|y|},\quad x,y\in \R,
\end{equation}where $C$ depends only on $f$.

The Bernstein space $B_\sigma$  is defined as the set  of entire functions $f$ satisfying (\ref{f}) with some $C$ depending only on $f$. An equivalent definition is that $B_\sigma$ consists of the bounded continuous functions that are the inverse Fourier transforms of tempered distributions supported by $[-\sigma,\sigma]$. 

Denote by $\mathcal{C}_0$ the collection of kernels $\{\varphi_u,u\in I\}$ satisfying the properties ($\beta$)-($\eta$) in the definition of $\mathcal{C}$ above. However, we do not require $I$ to be an interval. In particular, it can be a countable set.

\begin{theorem}\label{t2} Given a u.d. set $\La\subset\R$  and a kernel   $\{\varphi_u,u\in I\}\in\mathcal{C}_0$. The following conditions are equivalent:

{\rm(a)} For every $\sigma>0$ there is a constant $K=K(\sigma)$ such that
\begin{equation}\label{ssb}
\|f\|_\infty\leq K\sup_{\lambda\in\La, u\in I}|(f\ast\varphi_u)(\lambda)|,\quad \forall f\in B_\sigma;
\end{equation}

{\rm(b)}  $\La$ satisfies condition {\rm($\alpha$)}.
\end{theorem}

To prove this theorem we  need a lemma:

\begin{lemma}\label{l1}
Assume $f\in B_\sigma$ and $\{\varphi_u, u\in I\}\in \mathcal{C}_0$. If $(f\ast \varphi_u)(0)=0, u\in I,$ then  $f$ is odd, $f(-x)=-f(x),x\in\R$.
\end{lemma}

Proof. 1.  Given a function $f\in B_\sigma$, set
$$
f_r(z):=\frac{f(z)+\overline{f(\bar z)}}{2},\quad f_i(z):=\frac{f(z)-\overline{f(\bar z)}}{2i}.
$$ Then  $f_r,f_i$ are real (on $\R$) entire functions satisfying  $f=f_r+if_i$. It is clear that both $f_r$ and $f_i$ satisfy (\ref{f}), so that they  both lie in $B_\sigma$. Hence, since every $\varphi_u$ is real,  it suffices to prove the lemma for the real functions $f\in B_\sigma$.

2. Let us assume that  $f\in B_\sigma$ is real. Write
$$
f_e(x):=\frac{f(z)+f(- z)}{2},\quad f_o(x):=\frac{f(z)-f(- z)}{2}.
$$Clearly, $f_e\in B_\sigma$ is even, $f_o\in B_\sigma$ is odd and $f=f_e+f_o$. Since $\varphi_u$
is even, we have $(f_o\ast\varphi_u)(0)=0, u\in I$. Hence, to prove Lemma \ref{l1}, it suffices to check that if a real even function $f\in B_\sigma$ satisfies $(f\ast\varphi_u)(0)=0, u\in I$, then $f=0$.

3. Let us assume that $f\in B_\sigma$ is real, even and satisfies $(f\ast\varphi_u)(0)=0, u\in I$. If   $f$ does not vanish in $\Cc$ then $f(z)=e^{iaz}$ for some $-\sigma\leq a\leq\sigma,$ which implies $a=0, f(z)\equiv 1.$ Then $(f\ast\varphi_u)(0)=\Phi_u(0)=0, u\in I$, which contradicts condition ($\eta$).

Hence,  $f(w)=0$ for some $w\in\Cc.$ It follows that $f(-w)=0.$ Set $$g(z):=\frac{f(z)}{z^2-w^2}.$$Denote by $G$ the Fourier transform of $g$. Then $G$ is continuous, even and vanishes outside $(-\sigma,\sigma)$. Now, condition $(f\ast\varphi_u)=0,u\in I,$ implies:
$$
0=\int_\R\varphi_u(s)f(s)\,ds=\int_\R (s^2-w^2)\varphi_u(s)g(s)\,ds=$$$$-\int_{-\sigma}^\sigma (\Phi_u''(t)+w^2\Phi_u(t))G(t)\,dt=-2\int_{0}^\sigma (\Phi_u''(t)+w^2\Phi_u(t))G(t)\,dt.
$$
Using property ($\theta$), we conclude that $G=0$ and so $f=0$.

\subsection{Proof of Theorem \ref{t2}}

 We denote by $C$ different positive constants.

1. Suppose $W(\La)$ contains an empty set. It means that $\La$ contains arbitrarily long gaps: For every $\rho>0$ there exists $x_\rho$ such that
$\La\cap(x_\rho-2\rho,x_\rho+2\rho)=\emptyset$. Set \begin{equation}\label{ff}f_\rho(x):=\frac{\sin(\sigma(x-x_\rho))}{\sigma(x-x_\rho)}\in B_\sigma.\end{equation}Then  $\|f_\rho\|_\infty=1$. Using  (\ref{u1}), for all $x$ such that
$|x-x_\rho|\geq 2\rho,$ we have
$$
|(f_\rho\ast\varphi_u)(x)|\leq \int\limits_{|s|<\frac{|x-x_\rho|}{2}}\frac{2}{\sigma|x-x_\rho|}|\varphi_u(s)|\,ds+$$\begin{equation}\label{esti}\int\limits_{|s|>\frac{|x-x_\rho|}{2}}|\varphi_u(s)|\,ds\leq \frac{C}{|x-x_\rho|}.
\end{equation}
It readily follows that (\ref{ssb}) is not true.

2. Suppose $\La^\ast\subset b+(1/a)\Z$ for some $\La^\ast\in W(\La),b\in\R$ and $a>0$. Since $\La^\ast-b\in W(\La)$, we may assume that $b=0.$

Consider two cases: First, let us assume that
$\La\subset (1/a)\Z$. Set  $\sigma=\pi a$. Clearly,  the function $f(z):=\sin ( \pi a z)\in B_\sigma$. Since every function $\varphi_u$ is even while $f$ is odd, one may easily check  that
$(f\ast\varphi_u)(k/a)=0, k\in\Z$, so that  (\ref{ssb}) is not true.

Now, assume that $\La^\ast\subset (1/a)\Z$, for some $\La^\ast\in W(\La)$.
This means that for every small $\epsilon>0$  and large $R>0$ there is a point $v=v(\epsilon,R)\in \R$
such that $(\La-v)\cap(-R,R)$ is close to a subset of $(1/a)\Z$ in the sense that for every  $\lambda\in \La\cap(v-R,v+R)$ there exists $k(\lambda)\in\Z$ with  $$\label{la}|\lambda-v-k(\lambda)/a|\leq\epsilon,\quad \lambda\in \La\cap(v-R,v+R).$$

For simplicity of presentation, we assume that $v=0,a=1$,  and that\begin{equation}\label{la}\La\cap(-R,R)=\{\lambda_{k}: |k|\leq m\},\quad |\lambda_k-k|\leq\epsilon,\ \  m=[R], \ \  |k|\leq m.\end{equation}
The proof of the general case is similar.

Fix $\epsilon:=1/\sqrt R$. Set \begin{equation}\label{f1}f(x):=\sin ( \pi x)\frac{\sin(\epsilon x)}{\epsilon x}\in B_{\pi+\epsilon}\end{equation} and
$$
f_k(x):= \sin (\pi x)\frac{\sin(\epsilon \lambda_k)}{\epsilon \lambda_k}.
$$Then
$$
\left|f(\lambda_k-s)-(-1)^{k+1}f_k(s)\right|\leq\left|\left[\sin (\pi(\lambda_k-s))-\sin(\pi (k-s))\right]\frac{\sin\epsilon (\lambda_k-s)}{\epsilon (\lambda_k-s)}\right|+
$$\begin{equation}\label{00}
\left|\sin(\pi s)\left(\frac{\sin\epsilon (\lambda_k-s)}{\epsilon (\lambda_k-s)}-\frac{\sin\epsilon \lambda_k}{\epsilon \lambda_k}\right)\right|.
\end{equation}

By (\ref{la}),
$$
\left|\sin (\pi(\lambda_k-s))-\sin(\pi (k-s))\right|\leq \pi\epsilon,\quad s\in\R,
$$
 and so the first term in the right hand-side of (\ref{00}) is less than $\pi\epsilon$ for every $s\in\R$. To estimate the second term in (\ref{00}), we use the classical Bernstein's inequality (see e.g. \cite{ou1}, Lecture 2.10):
$$
\left|\frac{\sin\epsilon (\lambda_k-s)}{\epsilon (\lambda_k-s)}-\frac{\sin\epsilon \lambda_k}{\epsilon \lambda_k}\right|= \left|\int_0^s\left(\frac{\sin\epsilon (\lambda_k-u)}{\epsilon (\lambda_k-u)}\right)'du\right|\leq |s|
\left\|\left(\frac{\sin(\epsilon s)}{\epsilon s}\right)'\right\|_\infty\leq \epsilon |s|.
$$
Therefore, $$|f(\lambda_k-s)-(-1)^{k+1}f_k(s)|\leq \pi\epsilon(1+ |s|),\quad s\in\R.$$

Observe that
$$
(f\ast\varphi_u)(\lambda_k)=\int_\R(f(\lambda_k-s)-(-1)^{k+1}f_k(s))\varphi_u(s)\,ds+(-1)^{k+1}\int_\R f_k(s)\varphi_u(s)\,ds.
$$Since $f_k$ is odd, the last integral is equal to zero. It follows that for every $|k|\leq m$ we have
$$
|(f\ast\varphi_u)(\lambda_k)|\leq \pi\epsilon\int\limits_\R (1+|s|)|\varphi_u(s)|\,ds, \quad u\in I.$$
 Hence, using  (\ref{u1}) we conclude that $$|(f\ast\varphi_u)(\lambda)|\leq C\epsilon,\quad \lambda\in\La\cap(-R,R),\quad u\in I.$$

On the other hand, for all $\lambda\in\La, |\lambda|\geq R$ and $|s|<1/\epsilon =\sqrt R$, we get $$|f(\lambda-s)|\leq \frac{1}{\epsilon|\lambda-s|}\leq \frac{\sqrt R}{R-\sqrt R}<2\epsilon,$$provided $R$ is sufficiently large. This and (\ref{u1}) imply
$$|(f\ast\varphi_u)(\lambda)|\leq 2\epsilon\int\limits_{|s|<\sqrt R}|\varphi_u(s)|\,ds+\int\limits_{|s|>\sqrt{R}}|\varphi_u(s)|\,ds\leq C\epsilon,\quad \lambda\in\La,|\lambda|\geq R, u\in I.
$$
Since $\epsilon$ can be chosen arbitrarily small, we conclude that (\ref{ssb}) is not true.

3. Assume condition ($\alpha$) holds. We have to show that for every $\sigma>0$ there is a constant $K=K(\sigma)$ such that (\ref{ssb}) is true.

Assume this is not so. It means that there exists $\sigma>0$ and a sequence of functions $f_n\in B_\sigma$ satisfying
$$
\|f_n\|_\infty=1, \quad \sup_{u\in I,\lambda\in\La}|(f_n\ast\varphi_u)(\lambda)|\leq 1/n.
$$Choose points $x_n\in\R$ such that $|f_n(x_n)|>1-1/n$, and set $g_n(x):=f_n(x+x_n)$. It follows from the compactness property of Bernstein spaces (see e.g. \cite{ou1}, Lecture 2.8.3), that there is a subsequence $n_k$ such that  $g_{n_k}$ converge (uniformly on compacts in $\Cc$) to some non-zero function $g\in B_\sigma$. We may also assume (by taking if necessary a subsequence of $n_k$)  that the translates $\La-x_{n_k}$ converge weakly to some $\Gamma\in W(\La)$. By property ($\alpha$), $\Gamma$ is an  {\it infinite set} which is not a subset of any arithmetic progression.

Clearly, we have
$$
(g\ast\varphi_u)(\gamma)=0,\quad   u\in I, \ \gamma\in\Gamma.
$$By Lemma \ref{l1}, we see that every function $g(x-\gamma),\gamma\in
\Gamma,$ is odd. Clearly, this implies that $g$ is a periodic function and $\Gamma$ is a subset of  an arithmetic progression whose difference is a half-integer multiple of the period of $g$. Contradiction.

\section{Space--Time Sampling in Paley-Wiener Spaces}

Throughout this section we denote by $C$ different positive constants.

In what follows we assume that $I$  is an interval. We denote by $C$ different positive constants.

The following statement easily follows from (\ref{u2}) and
(\ref{ssb})  :
 \begin{corollary}\label{c1}
 Assume condition {\rm (\ref{ssb})} holds for some kernel $\{\varphi_u\}$ satisfying {\rm (\ref{u2})},   a u.d. set $\La$ and $\sigma>0$.
        Then there is a constant $K'=K'(\sigma)$ such that $$\|f\|_\infty^2\leq K'\int\limits_I\sup_{\lambda\in\La}|(f\ast\varphi_u)(\lambda)|^2\,du,\quad \forall f\in B_\sigma.$$
 \end{corollary}

We skip the simple proof.


\subsection{ Proof of Proposition \ref{p1}}

 Take any function $f\in PW_\sigma$ and denote by $F$ its  Fourier transform.
It follows from (\ref{fb}) that $\|\Phi_u\|_\infty\leq C, u\in I.$
Hence,  the functions $F\cdot\Phi_u\in L^2(-\sigma,\sigma)$ and $$\|f\ast\varphi_u\|_2=\|F\cdot\Phi_u\|_2\leq \|\Phi_u\|_\infty\|F\|_2\leq C\|f\|_2.$$

Clearly, $f\ast\varphi_u\in PW_\sigma,$ for every $u$. Using  Bessel's inequality (see e.g. Proposition~2.7 in \cite{ou1}), we get
$$
\sum_{\lambda\in\La}|(f\ast\varphi_u)(\lambda)|^2\leq C\|f\ast \varphi_u\|_2^2\leq C\|f\|_2^2,\quad u\in I,
$$which proves Proposition \ref{p1}.

\subsection{ Connection between space--time sampling in $B_{\sigma}$ and $PW_{\sigma}$}
Observe that if $\La$ is a sampling set (in the `classical sense') for the  Paley-Wiener space $PW_{\sigma'}$, then it is a sampling set for the  Bernstein spaces $B_{\sigma}$ with a `smaller' spectrum $\sigma<\sigma',$ and vice versa (see Theorem~3.32 in \cite{ou1}).
We provide a corresponding statement for the space-time sampling problem.

 For the reader's convenience, we recall the main inequalities:
\begin{equation}\label{l_ss_pw}
\|f\|^2_2\leq D \int\limits_I\sum_{\lambda\in\La}|(f\ast\varphi_u)(\lambda)|^2\,du,
\end{equation}
\begin{equation}\label{l_ss_b}
\|f\|_\infty\leq K\sup_{\lambda\in\La, u\in I}|(f\ast\varphi_u)(\lambda)|.
\end{equation}

\begin{theorem}\label{t3}
Let $\La$ be a u.d. set, a kernel $\{\varphi_u\}$ satisfy {\rm (\ref{u1})} and {\rm (\ref{u2})} and   $\sigma'>\sigma>0$.

{\rm (i)} Assume that {\rm (\ref{l_ss_b})} holds with some constant $K$ for all $f\in B_{\sigma'}$. Then there is a constant $D$ such that {\rm (\ref{l_ss_pw})} is true for every $f \in PW_{\sigma}$.

{\rm (ii)} Assume that {\rm (\ref{l_ss_pw})} holds with some constant $D$  for all $f \in PW_{\sigma'}$. Then there is a constant $K$ such that  {\rm (\ref{l_ss_b})} is true for every $f \in B_{\sigma}$.

\end{theorem}

{\noindent \bf Proof.}
The proof is somewhat similar to the proof of Theorem~3.32 in \cite{ou1}, but is more technical.

(i)  Assume that  (\ref{l_ss_b}) holds for every $f\in B_{\sigma'}$. Fix any positive number $\varepsilon$ satisfying \begin{equation}\label{epsi}\sigma+\varepsilon\leq\sigma'.\end{equation}

 Set
 \begin{equation}\label{h_def}
 \h(x):=\frac{\sin \varepsilon x}{\epsilon x},\quad \varepsilon>0.
 \end{equation}
It is easy to check that
 \begin{equation}\label{h}
 \h(0)=1,\quad \|\h\|_2^2=\frac{C}{\varepsilon},\quad \|\h'\|_2^2=C\varepsilon.
 \end{equation}

 For every $f\in PW_\sigma$, we have
 $$
 \|f\|_2^2=\int\limits_\R|f(x)|^2\,dx\leq \int\limits_\R\sup_{s\in\R}|\h(x-s)f(s)|^2\,dx.
 $$
Note that $\h(x-s)f(s) \in PW_{\sigma+\varepsilon}\subset B_{\sigma'}$. By Corollary 1, for every $x$ and $s$,   $$|\h(x-s)f(s)|^2\leq C\int\limits_I \sup_{\lambda\in\La}\left|\int\limits_\R\varphi_u(\lambda-s)\h(x-s)f(s)\,ds\right|^2du\leq $$
 $$ C\int\limits_I\sum_{\lambda\in\La}\left|\int\limits_\R\varphi_u(\lambda-s)\h(x-s)f(s)\,ds\right|^2du.$$

 Write$$J=J_u(x,\lambda):=\left|\int\limits_\R\varphi_u(\lambda-s)\h(x-s)f(s)\,ds\right|^2.
 $$Then
 \begin{equation}\label{fie}
 \|f\|^2_2\leq C\int\limits_\R \sum_{\lambda\in\La}\int\limits_I J\,dudx.  \end{equation}

 Clearly, $$J\leq 2(J_1+J_2),$$where $$J_1:=\left|\int\limits_\R\varphi_u(\lambda-s)\h(x-\lambda)f(s)\,ds\right|^2=|\h(x-\lambda)|^2\left|(f\ast\varphi_u)(\lambda)\right|^2,$$and using property (\ref{u1}) and the Cauchy–Schwartz inequality, we have
 $$
 J_2:=\left|\int\limits_\R\varphi_u(\lambda-s)(\h(x-s)-\h(x-\lambda))f(s)\,ds\right|^2\leq $$$$ \int\limits_\R|\varphi_u(s-\lambda)|\,ds\int\limits_\R|\varphi_u(\lambda-s)||\h(x-s)-\h(x-\lambda)|^2|f(s)|^2\,ds\leq$$
 $$C\int\limits_\R|\varphi_u(\lambda-s)||\h(x-s)-\h(x-\lambda)|^2|f(s)|^2\,ds.
 $$

  Observe that $$|\h(x-s)-\h(x-\lambda)|^2=\left|\int\limits_s^\lambda \h'(x-v)\,dv\right|^2\leq |s-\lambda|\int\limits_s^\lambda|\h'(x-v)|^2\,dv.$$Hence,
 $$
 J_2\leq C \int\limits_\R|\varphi_u(\lambda-s)||s-\lambda|\left(\int\limits_s^\lambda|\h'(x-v))|^2dv\right) |f(s)|^2\,ds.
 $$

 Using (\ref{h}), we have
 $$
 \int\limits_\R\sum_{\lambda\in\La}\int\limits_I J_1\,dudx=\int\limits_\R|\h(\lambda-x)|^2\,dx\sum_{\lambda\in\La}\int\limits_I |(f\ast\varphi_u)(\lambda)|^2du\leq $$$$ \frac{C}{\epsilon}\sum_{\lambda\in\La}\int\limits_I |(f\ast\varphi_u)(\lambda)|^2du.
 $$

 To estimate the second sum we switch the order of integration and apply (\ref{h}):
 $$
 \int\limits_\R\sum_{\lambda\in\La} \int\limits_I J_2\,dudx\leq $$$$\int\limits_\R \sum_{\lambda\in\La}\int\limits_I|\varphi_u(\lambda-s)||s-\lambda||f(s)|^2\left(\int\limits_\R\int\limits_s^\lambda|\h'(x-v)|^2dv\,dx\right)duds\leq$$$$ C\varepsilon \int\limits_\R\int\limits_I\sum_{\lambda\in\La} |\varphi_u(\lambda-s)||s-\lambda|^2|f(s)|^2duds.
 $$

Now,  by  (\ref{u1}) we get$$
\sum_{\lambda\in\La} |\varphi_u(\lambda-s)||s-\lambda|^2\leq C\sum_{\lambda\in\La} \frac{(\lambda-s)^2}{1+(\lambda-s)^4}<C,\quad u\in I,s\in\R.
$$where the second inequality holds since $\La$ is a u.d. set (see definition in (\ref{sep})).  Hence,$$\int\limits_\R\sum_{\lambda\in\La} \int\limits_I J_2\,dudx\leq C\varepsilon |I|\|f\|_2^2,$$where $|I|$ is the length of $I$.

 Combining  this with the  estimate for $J_1$ and using (\ref{fie}), we conclude that
 $$
 \|f\|_2^2\leq \frac{C}{\varepsilon}\sum_{\lambda\in\La}\int\limits_I|(f\ast\varphi_u)(\lambda)|^2\,du+C\varepsilon |I|\|f\|_2^2.
 $$
 Choosing $\varepsilon$ small enough, we obtain (\ref{l_ss_pw}).

(ii)  Assume  (\ref{l_ss_pw}) holds with some constant $D$  for all $f \in PW_{\sigma'}$.

 We will argue by contradiction. Assume that  there is no constant $K$ such that (\ref{l_ss_b}) holds for every $f\in B_{\sigma}$. This means that there exist $g_j \in B_{\sigma}$ such that $\|g_j\|_\infty=1$,
\begin{equation}\label{contr_eq}
    \sup\limits_{u \in I, \lambda \in \La} \left|(g_j\ast \varphi_u) (\lambda)\right| < \frac{1}{j},
\end{equation}
and for some points $x_j$ we have $|g_j(x_j)|\geq 1/2$.

Assume $\varepsilon>0$ satisfies (\ref{epsi}) and let $h_\varepsilon$ be defined by formula~(\ref{h_def}). Set $$f_j(x): = g_j(x) h_{\varepsilon}(x - x_j).$$ It is clear that for every $j$ we have $f_j\in PW_{\sigma'}$, $\|f_j\|_\infty\leq 1$, and that $|f_j(x_j)|\geq 1/2$.
The last two inequalities and the Bernstein's inequality  imply that there is a constant  $K'>0$ such that  \begin{equation}\label{berns} \|f_j\|_2 \ge K',\quad j\in\N.\end{equation}

By (\ref{l_ss_pw}), we get
$$
\|f_j\|^2_{2} \le C \int\limits_I \sum\limits_{\lambda \in \La} |(f_j\ast\varphi_u)(\lambda)|^2 du = C\int\limits_I \sum\limits_{\lambda \in \La} \left|\int_\R g_j(x)\varphi_u(\lambda-x)\h(x-x_j) dx \right|^2 du.
$$This gives
\begin{equation}\label{sist}
\|f_j\|_2^2\leq C (\J_1+\J_2),
\end{equation}where $\J_1$ and $\J_2$ are defined as follows:
$$
\J_1: = \int\limits_I \sum\limits_{\lambda \in \La} \left|\int\limits_\R g_j(x)\varphi_u(\lambda-x)(\h(x-x_j) - \h(\lambda - x_j) ) dx\right|^2\, du,
$$
$$
\J_2 :=  \int\limits_I \sum\limits_{\lambda \in \La} \left|\int\limits_\R g_j(x)\varphi_u(\lambda-x) \h(\lambda - x_j) dx\right|^2 \, du.
$$

By Bessel's inequality (see, e.g. \cite{ou1}, Proposition 2.7) and (\ref{h}),
$$
\sum\limits_{\lambda \in \La} |\h (\lambda - s)|^2 \le C\|\h\|_2^2\leq \frac{C}{\varepsilon},\quad \forall s\in\R.
$$Therefore, using (\ref{contr_eq}) we arrive at
$$\J_2 \le  \frac{C}{\varepsilon j^2}|I|.$$

Let us now estimate $\J_1$. Recall that  $\|g_j\|_{\infty} = 1$. Using the change of variables $x=t+\lambda$,
we get
$$
\J_1\leq \int\limits_I\sum_{\lambda\in\La}\left(\int\limits_\R\left|\varphi_u(-t)\int\limits_0^t \h'(s+\lambda-x_j)\,ds\right|dt\right)^2\,du.
$$
Now, use  the Cauchy–Schwarz inequality:
$$
\J_1\leq \int\limits_I\sum_{\lambda\in\La}\int\limits_\R|\varphi_u(-t)|^2(1+t^2)^2\,dt\int\limits_\R\frac{1}{(1+t^2)^2}\left|\int\limits_0^t \h'(s+\lambda-x_j)\,ds\right|^2dt\,du.
$$Using again the Cauchy–Schwarz inequality and condition (\ref{u1}), we arrive at
 $$
\J_1\leq C\int\limits_I\sum_{\lambda\in\La}\int\limits_\R\frac{|t|}{(1+t^2)^2}\left| \int\limits_0^t |\h'(s+\lambda-x_j)|^2\,ds\right|\,dt\,du.
$$Finally, Bessel's inequality yields
$$
\sum_{\lambda\in\La}|\h'(s+\lambda-x_j)|^2\leq C\|\h\|_2^2\leq C\varepsilon,
$$and we conclude that $$\J_1\leq C|I|\varepsilon.$$

 We now insert the estimate for $\J_1,\J_2$ in (\ref{sist}) and use (\ref{berns}) to get the estimate
$$
(K')^2 \le \frac{C}{\varepsilon j^2} + C|I| \varepsilon.
$$ Choosing  $\varepsilon$  sufficiently small, we arrive at contradiction for all large enough  $j$.

\subsection{ Proof of Theorem \ref{t1}}
The proof easily follows from Theorems \ref{t2} and \ref{t3}.


Assume that the assumptions of Theorem \ref{t1} hold.

(i) Assume that $\Lambda$ satisfies condition $(\alpha)$. Then by Theorem~{\ref{t2}}, for every $\sigma>0$ there exists $K=K(\sigma)$ such that inequality (\ref{ssb}) is true. Applying Theorem \ref{t3}, we see that there exists $A=A(\sigma)>0$ the left hand-side inequality in (\ref{ss}) is also true for every $\sigma>0$.


(ii) Assume that $\Lambda$ does not satisfy condition  $(\alpha)$. Then by  Theorem~\ref{t2}, there exists $\sigma>0$ such that there is no constant $K$ for which condition (\ref{ssb}) is true.
Applying Theorem~\ref{t3}, we see that for every positive $\sigma'>\sigma$ there is no constant $D$ such that inequality
(\ref{l_ss_pw}) holds for every $f\in PW_{\sigma'}$.

\noindent Alexander Ulanovskii\\ University of Stavanger, Department of Mathematics and Physics,\\
4036 Stavanger, Norway,\\
alexander.ulanovskii@uis.no

\medskip

\noindent Ilya Zlotnikov\\
University of Stavanger, Department of Mathematics and Physics,\\
4036 Stavanger, Norway,\\
ilia.k.zlotnikov@uis.no

\end{document}